\documentclass[12pt]{article}
\usepackage{caption}
\usepackage{mathrsfs}
\usepackage{amsmath,amsfonts,amssymb,rotating,amsthm}
\usepackage{bm}
\usepackage{hyperref}
\usepackage{array}
\usepackage{color}
\usepackage{breqn}
\usepackage{enumitem}
\usepackage{listings} 
\usepackage{xcolor} 
\usepackage{booktabs} 
\usepackage{multirow}
\usepackage{tikz-cd}
\usepackage{tikz}
\allowdisplaybreaks[4]
\newtheorem{theo}{Theorem}[section]
\newtheorem{deff}[theo]{Definition}

\theoremstyle{remark}

\numberwithin{equation}{section}

\newdimen\Squaresize \Squaresize=11pt
\newdimen\Thickness \Thickness=0.7pt
\def\Square#1{\hbox{\vrule width \Thickness
		\vbox to \Squaresize{\hrule height \Thickness\vss
			\hbox to \Squaresize{\hss#1\hss}
			\vss\hrule height\Thickness}
		\unskip\vrule width \Thickness} \kern-\Thickness}

\def\Vsquare#1{\vbox{\Square{$#1$}}\kern-\Thickness}

\def\moins{\raise 1pt\hbox{{$\scriptstyle -$}}}

\newcommand{\x}{\ensuremath{{\bm x}}}
\newcommand{\y}{\ensuremath{{\bm y}}}

\begin{document}

	\begin{center}
	{\Large \bf New lower bounds on the non-repetitive chromatic number of some graphs}
\end{center}

\begin{center}
	Tianyi Tao$^{1}$, Junchi Zhang$^{2}$, Wentao Zhang$^{3}$, and Alex Toole$^4$\\[8pt]
	$^{1,2,3}$School of Mathematical Sciences\\
	Fudan University\\
	Shanghai, China\\[6pt]
	$^4$Mathematics Department\\
	Palomar College\\
	 San Marcos, California, USA\\[6pt]
	Email: {\tt $^{1}$tytao20@fudan.edu.cn, $^{2}$jczhang24@m.fudan.edu.cn,  $^{3}$wtzhang20@fudan.edu.cn, $^4$lxtoole@gmail.com\tt}
\end{center}

\vspace{0.3cm} \noindent{\bf Abstract.} 
A graph \( G \) is said to be (vertex) non-repetitively colored if no simple path in \( G \) has a sequence of vertex colors that forms a repetition. Formally, a coloring \( c: V(G) \to \{1, 2, \dots, k\} \) is non-repetitive if, for every path \(\langle v_1, v_2, \dots, v_{2m} \rangle\) in \( G \), the sequence of colors \( c(v_1), c(v_2), \dots, c(v_{2m}) \) is not of the form \( ww \), where \( w \) is a sequence of \( m \) colors. The minimum number of colors required for such a coloring is called the \emph{non-repetitive chromatic number} of \(G\), denoted by \(\pi(G)\). In this paper, we primarily prove that \(\pi(P \square P) \ge 6\) and \(\pi(P \boxtimes P) \ge 9\), where \( P \square P \) and \( P \boxtimes P \) are the Cartesian product and the strong product of two infinite paths, respectively. This improves upon the previous best lower bounds.
\\

\noindent {\bf Keywords:} non-repetitive coloring; algorithm; grid graph.
\\
\noindent {\bf 2020 Mathematics Subject Classification:} 05C85, 05C15, 05C76.
\\

\section{Introduction}

In this paper, we denote by \( P_n \) the path with \( n \) vertices, and by \( P \) the infinite path. Specifically, \( V(P) = \mathbb{Z} \), and \( i \sim j \) if and only if \( |i - j| = 1 \).

\begin{deff}\cite{alon2002nonrepetitive}
	Let \( G \) be a graph. Suppose \( f \) is a vertex coloring of \( G \), \( k \) is a positive integer, and \( P_{2k} = \langle v_1, \dots, v_k, v_{k+1}, \dots, v_{2k} \rangle \) is a path of \( 2k \) vertices in \( G \). We say that \( f \) \textit{repetitively colors} \( P_{2k} \) if for all \( 1 \leq i \leq k \), \( f(v_i) = f(v_{k+i}) \). If no even-length path in \( G \) is repetitively colored by \( f \), then \( f \) is called a \textit{non-repetitive coloring} of \( G \). The minimum number of colors required for \( G \) to admit such a coloring is denoted by \( \pi(G) \), referred to as the \textit{non-repetitive chromatic number} or \textit{Thue number} of \( G \).
	
	An even-length path is called a repetitive path if it is repetitively colored by some $f$.
\end{deff}

The study of problems in this field originates from the following theorem.

\begin{theo}[Thue~\cite{thue1906uber}]
	\( \pi(P) = 3 \). More specifically, for all \( n \geq 4 \), \( \pi(P_n) = 3 \).
\end{theo}

For additional literature on non-repetitive coloring, we can refer to \cite{barat2005NoNG,barat2007square,BOTLER2025221,brevsar2007nonrepetitive,dujmovic2020planar,grytczuk2007nonrepetitive,harant2012nonrepetitive,KESZEGH2020108,WoodDavidR2021NGC}.

Dujmović, Joret, Micek, Morin, Ueckerdt, and Wood\cite{dujmovicVida2020Pghb} proved that every planar graph \( G \) is a subgraph of \( H \boxtimes P \boxtimes K_3 \) for some graph \( H \) with treewidth at most 3 and some path \( P \). This result implies that \( \pi(G) \leq 768 \) (see also \cite{kundgen2008nonrepetitive}). The proof of this theorem motivated our study of non-repetitive colorings in graph products.

\begin{deff}
	The \textit{Cartesian product} of two graphs \( G \) and \( H \), denoted by \( G \square H \), is a graph whose vertex set is the Cartesian product \( V(G) \times V(H) \). Two vertices \( (u_1, v_1) \) and \( (u_2, v_2) \) are adjacent in \( G \square H \) if and only if either:
	\begin{itemize}
		\item \( u_1 = u_2 \) and \( v_1 \) is adjacent to \( v_2 \) in \( H \), or
		\item \( v_1 = v_2 \) and \( u_1 \) is adjacent to \( u_2 \) in \( G \).
	\end{itemize}
\end{deff}

\begin{deff}
	The \textit{strong product} of two graphs \( G \) and \( H \), denoted by \( G \boxtimes H \), is a graph whose vertex set is the Cartesian product \( V(G) \times V(H) \). Two vertices \( (u_1, v_1) \) and \( (u_2, v_2) \) are adjacent in \( G \boxtimes H \) if and only if one of the following holds:
	\begin{itemize}
		\item \( u_1 = u_2 \) and \( v_1 \) is adjacent to \( v_2 \) in \( H \),
		\item \( v_1 = v_2 \) and \( u_1 \) is adjacent to \( u_2 \) in \( G \), or
		\item \( u_1 \) is adjacent to \( u_2 \) in \( G \) and \( v_1 \) is adjacent to \( v_2 \) in \( H \).
	\end{itemize}
\end{deff}

Now, we provide a simple estimate for the non-repetitive chromatic numbers of \( P \square P \) and \( P \boxtimes P \). Clearly, \( \pi(P \boxtimes P) \ge \pi(P \square P) \). Since \( P_4 \) is a subgraph of \( P \square P \), we have \( \pi(P \boxtimes P) \ge \pi(P \square P) \ge \pi(P_4) = 3 \). Moreover, because \( C_{10} \) is a subgraph of \( P \square P \), the following theorem implies \( \pi(P \square P) \ge 4 \):

\begin{theo}[Currie~\cite{currie2002there}]
	For an \( n \)-vertex cycle \( C_n \),
	\[
	\pi(C_n) = 
	\begin{cases} 
		4 & \text{if } n \in \{5, 7, 9, 10, 14, 17\}, \\ 
		3 & \text{otherwise}.
	\end{cases}
	\]
\end{theo}

Furthermore, regarding the non-repetitive chromatic number of \( P \square P \), there exist the following non-trivial upper and lower bounds.

\begin{theo}[Tao~\cite{tao2024nonrepetitivecoloringsgrids}]
	\(5\le\pi(P\square P)\le 12.\)
\end{theo}

It is worth mentioning that Fertin, Raspaud, and Reed \cite{FertinGuillaume2004Scog} claimed to have computationally verified that the star chromatic number of \( P \square P \) is 5, which also implies \( \pi(P \square P) \geq 5 \). This also follows from a case analysis in~\cite{toole} that establishes \( \pi(P_2 \square P_n) \geq 5 \) for \( n\ge 9 \), and provides a coloring that suggests that equality holds. 

About the upper bound for \( \pi(P \boxtimes P) \), the best known result is given by the following theorem:

\begin{theo}[K\"undgen, Pelsmajer~\cite{kundgen2008nonrepetitive}]
	\( \pi(P \boxtimes P)\le 16. \)
\end{theo}

In this paper, we will design an algorithm to provide improved lower bounds for \( \pi(P \square P) \) and \( \pi(P \boxtimes P) \). Specifically, we prove the following theorem, which gives a positive answer to Question 5.2 in \cite{tao2024nonrepetitivecoloringsgrids}:

\begin{theo}\label{1}
	\( \pi(P \square P)\ge 6 \) and \( \pi(P \boxtimes P)\ge 9.\)
\end{theo}

Moreover, we will also study a graph that lies between \( P\square P \) and \( P \boxtimes P \). The triangular tiling graph \( T_3 \) is defined as follows:  
\( T_3 \) is an infinite graph where each vertex has degree 6, and the edges form a regular tiling of the Euclidean plane with equilateral triangles. Each vertex is adjacent to six other vertices, forming a hexagonal neighborhood structure. Figure 1 is a schematic diagram of $T_3$.

\begin{figure}[htbp]
	\centering
	\begin{tikzpicture}[
		tri/.style={thick, draw=black},
		dot/.style={circle, fill=black, inner sep=0pt, minimum size=4pt}
		]
		
		\def\side{1}
		\pgfmathsetmacro{\hgt}{sqrt(3)/2*\side}
		
		\foreach \x in {0,1,...,4} {
			\foreach \y in {0,2,4} {
				\node[dot] at (\x,\hgt*\y) {};
				\draw[tri] (\x,\hgt*\y) -- (\x-1/2,\hgt*\y+\hgt);
			}
		}
		
		\foreach \x in {0,1,...,3} {
			\foreach \y in {0,2,4} {
				\draw[tri] (\x,\hgt*\y) -- (\x+1,\hgt*\y);
				\draw[tri] (\x,\hgt*\y) -- (\x+1/2,\hgt*\y+\hgt);
			}
		}
		
		\foreach \x in {0,1,...,4} {
			\foreach \y in {1,3} {
				\node[dot] at (\x-1/2,\hgt*\y) {};
			}
		}
		
		\foreach \x in {0,1,...,3} {
			\foreach \y in {1,3} {
				\draw[tri] (\x-1/2,\hgt*\y) -- (\x+1/2,\hgt*\y);
			}
		}
		
		\foreach \x in {0,1,...,4} {
			\foreach \y in {1,3} {
				\draw[tri] (\x-1/2,\hgt*\y) -- (\x,\hgt*\y+\hgt);
			}
		}
		
		\foreach \x in {0,1,...,3} {
			\foreach \y in {1,3} {
				\draw[tri] (\x+1/2,\hgt*\y) -- (\x,\hgt*\y+\hgt);
			}
		}
		
		\foreach \x in {0,1,...,4} {
			\node[dot] at (\x-1/2,\hgt*5) {};
		}
		
		\foreach \x in {0,1,...,3} {
			\draw[tri] (\x-1/2,\hgt*5) -- (\x+1/2,\hgt*5);
		}
		
	\end{tikzpicture}
	\caption{}
	\label{fig:triangular_lattice}
\end{figure}

It can be observed that $P \square P$ is contained in \( T_3 \), which in turn is contained in \( P \boxtimes P \). We provide a lower bound for the non-repetitive chromatic numbers of \( T_3 \) :

\begin{theo}\label{2}
	\( \pi(T_3)\ge 9 \).
\end{theo}

\section{Algorithm}\label{sec-2}

The problem of determining whether a given coloring of a given graph is non-repetitive was addressed by Toole~\cite{toole}, who provided a specific algorithm. We describe it below:

\begin{enumerate}
	\item Input the graph \( G \) and its coloring \( f \). For each pair of vertices \( (u, v) \) with the same color, initialize \( u_1 = u \), \( v_1 = v \), and set \( {\bf P}_1 = \langle u_1 \rangle \), \( {\bf P}_2 = \langle v_1 \rangle \). For \( k \geq 1 \), given \( {\bf P}_1 = \langle u_1, \dots, u_k \rangle \) and \( {\bf P}_2 = \langle v_1, \dots, v_k \rangle \), check if \( v_k \) is adjacent to \( u_1 \) or \( u_k \) is adjacent to \( v_1 \). If so, a repetitive path is detected, and the algorithm outputs False.
	
	If neither \( v_k \) is adjacent to \( u_1 \) nor \( u_k \) is adjacent to \( v_1 \), search for vertices \( u_{k+1} \) and \( v_{k+1} \) in the neighborhoods of \( u_k \) and \( v_k \) respectively, such that:
	\begin{itemize}
		\item \( u_{k+1} \) and \( v_{k+1} \) have the same color,
		\item \( u_{k+1} \neq v_{k+1} \),
		\item \( \{u_{k+1}, v_{k+1}\} \cap (V({\bf P}_1) \cup V({\bf P}_2)) = \varnothing \).
	\end{itemize}
	If such \( u_{k+1} \) and \( v_{k+1} \) exist, update \( {\bf P}_1 \) to \( \langle u_1, \dots, u_k, u_{k+1} \rangle \) and \( {\bf P}_2 \) to \( \langle v_1, \dots, v_k, v_{k+1} \rangle \), and continue the search. If no such vertices exist, remove \( u_k \) from \( {\bf P}_1 \) and \( v_k \) from \( {\bf P}_2 \) (i.e., backtrack to \( {\bf P}_1 = \langle u_1, \dots, u_{k-1} \rangle \) and \( {\bf P}_2 = \langle v_1, \dots, v_{k-1} \rangle \)).
	
	We stress that the paths ${\bf P}_1$ and ${\bf P}_2$ have identical sequences of colors.
	
	\item Iterate over all pairs of vertices \( (u, v) \) with the same color. If the algorithm completes without outputting False, then \( f \) is a non-repetitive coloring, and the algorithm outputs True.
\end{enumerate}

Note that determining whether a given coloring of a graph is \textnormal{non-repetitive} is a \textbf{coNP-complete} problem, even when only 4 colors are used. For details, see \cite{MARX200913}.

For a finite or infinite graph \( G \), we now aim to find a lower bound for the non-repetitive chromatic number of \( G \). To achieve this, we need to verify that for some \( c \), \( G \) does not admit a non-repetitive \( c \)-coloring, which implies \( \pi(G) \geq c + 1 \). Clearly, due to limited computational power, it is impossible to enumerate all colorings of \( G \) and demonstrate that every \( c \)-coloring contains a repetitive path. Therefore, we design the following algorithm to compute the number of non-repetitive \( c \)-colorings for certain subgraphs \( H \) of \( G \), and gradually increase the size of \( H \) until the number of non-repetitive \( c \)-colorings becomes zero.

\begin{enumerate}
	\item Input the graph \( G \) and a positive integer \( c \). Our goal is to prove that \( \pi(G) \geq c + 1 \). Select a suitable subgraph \( H \) of \( G \), store all non-repetitive \( c \)-colorings of \( H \), and output the number of non-repetitive \( c \)-colorings of \( H \).
	
	\item Select a vertex \( v \in V(G) \setminus V(H) \). Iterate over all non-repetitive \( c \)-colorings of \( H \). For each coloring, assign \( v \) each of the \( c \) possible colors, and check whether the induced subgraph \( G[V(H) \cup \{v\}] \) contains a repetitive path. If no repetitive path is found, store this coloring and update \( H \) to \( G[V(H) \cup \{v\}] \). Output the number of non-repetitive \( c \)-colorings of the updated \( H \).
	
	\item Repeat the process by selecting suitable vertices \( v \) until a subgraph \( H \) of \( G \) is obtained for which the number of non-repetitive \( c \)-colorings is zero. The algorithm then terminates.
\end{enumerate}

Under this framework, the search space is significantly reduced. In this algorithm, if we choose an appropriate initial subgraph \( H \) and a suitable order for adding vertices, the algorithm can run much faster. In the following content, we will describe how to determine \( H \) and the order of adding vertices \( v \) for specific graphs.

Additionally, note that when checking whether \( G[V(H) \cup \{v\}] \) contains a repetitive path, we only need to verify if there exists a repetitive path that includes \( v \), since \( H \) is already non-repetitively colored. This observation greatly reduces the computational cost. To address this, we need to modify Toole's algorithm accordingly:

Let \( v \in V(G) \), and let \( f \) be a \( c \)-coloring of \( V(G) \). Our goal is to determine whether there exists no repetitive path containing \( v \).

\begin{enumerate}
	\item Input a graph \( G \), a coloring \( f \) of \( V(G) \), and a specified vertex \( v \in V(G) \). For each pair of vertices \( (u, v) \) with the same color, initialize two paths \( {\bf P}_1 = \langle u \rangle \) and \( {\bf P}_2 = \langle v \rangle \). For \( {\bf P}_1 = \langle x_1, \dots, x_k \rangle \) and \( {\bf P}_2 = \langle y_1, \dots, y_k \rangle \), if \( x_k \) is adjacent to \( y_1 \) or \( x_1 \) is adjacent to \( y_k \), then a repetitive path is detected, and the algorithm outputs False.
	
	Otherwise, iteratively check all possible path extensions by selecting neighbors $x$ of $x_1$ and $y$ of $y_1$, or neighbors $x$ of $x_k$ and $y$ of $y_k$, such that:
	\begin{itemize}
		\item $x$ and $y$ have the same color,
		\item $x \neq y$,
		\item $\{x, y\} \cap (V({\bf P}_1) \cup V({\bf P}_2)) = \varnothing$.
	\end{itemize}
	If such $x$ and $y$ exist, update ${\bf P}_1$ and ${\bf P}_2$ as follows:
	\begin{itemize}
		\item If $x$ is a neighbor of $x_1$ and $y$ is a neighbor of $y_1$, update ${\bf P}_1$ to $\langle x, {\bf P}_1 \rangle$ and ${\bf P}_2$ to $\langle y, {\bf P}_2 \rangle$.
		\item If $x$ is a neighbor of $x_k$ and $y$ is a neighbor of $y_k$, update ${\bf P}_1$ to $\langle {\bf P}_1, x \rangle$ and ${\bf P}_2$ to $\langle {\bf P}_2, y \rangle$.
	\end{itemize}
	If no such vertices exist, remove $x_k$ from ${\bf P}_1$ and $y_k$ from ${\bf P}_2$ (i.e., backtrack to ${\bf P}_1 = \langle x_1, \dots, x_{k-1} \rangle$ and ${\bf P}_2 = \langle y_1, \dots, y_{k-1} \rangle$).
	
	\item Iterate over all vertices \( u \) with the same color as \( v \). If the algorithm completes without outputting False, then no repetitive path contains \( v \), and the algorithm outputs True, which means there is no repetitive path containing \( v \).
\end{enumerate}

Now we discuss how to select an appropriate initial subgraph \( H \) and the order of adding vertices for a specific graph \( G \). Taking \( G = P \square P \) and \( c = 5 \) as an example, we choose \( H = C_4 \). Under symmetry, there are only three distinct non-repetitive colorings of \( H \), as shown in Figure 2 (numbers represent colors).

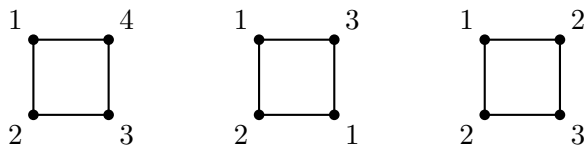
\begin{figure}[h]
	\centering
	\begin{tikzpicture}[
		dot/.style={circle, fill=black, inner sep=0pt, minimum size=4pt},
		label/.style={font=\small, text=black}
		]

		\begin{scope}[xshift=0cm] 
			\draw[thick] (0,0) rectangle (1cm,1cm);

			\node[dot] at (0,1) {};
			\node[dot] at (1,1) {};
			\node[dot] at (0,0) {};
			\node[dot] at (1,0) {};

			\draw (0,1) node [above left, font=\small] {1};
			\draw (0,0) node [below left, font=\small] {2};
			\draw (1,0) node [below right, font=\small] {3};
			\draw (1,1) node [above right, font=\small] {4};
		\end{scope}

		\begin{scope}[xshift=3cm] 
			\draw[thick] (0,0) rectangle (1cm,1cm);

			\node[dot] at (0,1) {};
			\node[dot] at (1,1) {};
			\node[dot] at (0,0) {};
			\node[dot] at (1,0) {};

			\draw (0,1) node [above left, font=\small] {1};
			\draw (0,0) node [below left, font=\small] {2};
			\draw (1,0) node [below right, font=\small] {1};
			\draw (1,1) node [above right, font=\small] {3};
		\end{scope}
		
		\begin{scope}[xshift=6cm] 
			\draw[thick] (0,0) rectangle (1cm,1cm);

			\node[dot] at (0,1) {};
			\node[dot] at (1,1) {};
			\node[dot] at (0,0) {};
			\node[dot] at (1,0) {};

			\draw (0,1) node [above left, font=\small] {1};
			\draw (0,0) node [below left, font=\small] {2};
			\draw (1,0) node [below right, font=\small] {3};
			\draw (1,1) node [above right, font=\small] {2};
		\end{scope}
		
	\end{tikzpicture}
	\caption{Non-repetitive colorings of \( H = C_4 \).}
	\label{fig:figure2}
\end{figure}

We denote by \( v_i \) the vertex added at the \( i \)-th iteration of the algorithm. It is easy to see that when \( v \) is a degree-1 vertex in \( G[V(H) \cup \{v\}] \), any repetitive path containing \( v \) can only start (or end) at \( v \). As the degree of \( v \) increases, the number of possible paths that need to be checked also increases. Therefore, we adopt the approach shown in Figure 3, where each added \( v_i \) is either a degree-1 or degree-2 vertex. The top-left four unlabeled vertices represent the initial subgraph \( H \).

\begin{figure}[h]
	\centering
	\begin{tikzpicture}[
		dot/.style={circle, fill=black, inner sep=0pt, minimum size=4pt},
		label/.style={font=\small},
		highlight/.style={circle, draw=red, thick, minimum size=3cm} 
		]

		\coordinate (A) at (0, 3);     
		\coordinate (B) at (1.5, 3);  
		\coordinate (C) at (3, 3);      
		\coordinate (D) at (4.5, 3);   
		\coordinate (E) at (6, 3);    
		
		\coordinate (F) at (0, 1.5);   
		\coordinate (G) at (1.5, 1.5);  
		\coordinate (H) at (3, 1.5);   
		\coordinate (I) at (4.5, 1.5);
		\coordinate (J) at (6, 1.5);    
		
		\coordinate (K) at (0, 0);      
		\coordinate (L) at (1.5, 0);   
		\coordinate (M) at (3, 0);     
		\coordinate (N) at (4.5, 0);   
		\coordinate (O) at (6, 0);     
		
		\coordinate (P) at (0, -1.5);  
		\coordinate (Q) at (1.5, -1.5); 
		\coordinate (R) at (3, -1.5);   
		\coordinate (S) at (4.5, -1.5); 
		\coordinate (T) at (6, -1.5);   
		
		\node[dot] at (A) {}; 
		\node[dot] at (B) {}; 
		\node[dot] at (C) {}; 
		\node[dot] at (D) {}; 
		\node[dot] at (E) {}; 
		\node[dot] at (F) {}; 
		\node[dot] at (G) {}; 
		\node[dot] at (H) {}; 
		\node[dot] at (I) {}; 
		\node[dot] at (K) {}; 
		\node[dot] at (L) {}; 
		\node[dot] at (M) {}; 
		\node[dot] at (N) {}; 
		\node[dot] at (P) {}; 
		\node[dot] at (Q) {}; 
		\node[dot] at (R) {}; 
		\node[dot] at (S) {}; 
		
		\draw (A) node[below left] {}; 
		\draw (B) node[below left] {}; 
		\draw (C) node[below left] {$v_1$};
		\draw (D) node[below left] {$v_6$};
		\draw (E) node[below left] {$v_{13}$};
		\draw (F) node[below left] {};
		\draw (G) node[below left] {};
		\draw (H) node[below left] {$v_2$};
		\draw (I) node[below left] {$v_7$};
		\draw (J) node[below left] {$\cdots$};
		\draw (K) node[below left] {$v_5$};
		\draw (L) node[below left] {$v_4$};
		\draw (M) node[below left] {$v_3$};
		\draw (N) node[below left] {$v_8$};
		\draw (P) node[below left] {$v_{12}$};
		\draw (Q) node[below left] {$v_{11}$};
		\draw (R) node[below left] {$v_{10}$};
		\draw (S) node[below left] {$v_9$};
		
		\draw (A) -- (B); 
		\draw (B) -- (C); 
		\draw (C) -- (D); 
		\draw (D) -- (E); 
		
		\draw (A) -- (F); 
		\draw (B) -- (G); 
		\draw (C) -- (H); 
		\draw (D) -- (I);

		\draw (F) -- (G); 
		\draw (G) -- (H); 
		\draw (H) -- (I);

		\draw (F) -- (K); 
		\draw (G) -- (L); 
		\draw (H) -- (M);
		\draw (I) -- (N); 
		
		\draw (K) -- (L); 
		\draw (L) -- (M); 
		\draw (M) -- (N);

		\draw (K) -- (P); 
		\draw (L) -- (Q); 
		\draw (M) -- (R); 
		\draw (N) -- (S); 
		
		\draw (P) -- (Q); 
		\draw (Q) -- (R); 
		\draw (R) -- (S);

		\node[highlight] at (0.75, 2.25) {}; 
		\draw (0.75, 2.25) node[above right, red] {$H$}; 
	\end{tikzpicture}
	\caption{}
\end{figure}

Note that the output of each step of the algorithm is not the exact number of non-repetitive colorings of \( H \), but rather the number up to certain symmetries. This is because when assigning colors to the vertices of the initial \( H \), we treat colorings that are equivalent under permutation symmetry as the same. However, this perspective does not affect our results and, in fact, reduces some computational complexity.

Let \( n(i) \), \( i \geq 4 \), print the number of non-repetitive $c$-colorings of \( H \) when \( |H| = i \). In Section 3, we will present computational results for this function.

For \( P \boxtimes P \) and \( T_3 \), note that \( T_3 \) is a subgraph of \( P \boxtimes P \). We adopt the same approach as for \( P \square P \), starting with an initial subgraph with 4 vertices and setting the vertices \( v_i \) in the same manner as in Figure 3, increasing the value of \( c \) until the best lower bound provided by the algorithm is obtained. These computational results will be presented in Section 3.

\section{Results}

\begin{table}[htbp]
	\centering
	\begin{minipage}{0.48\textwidth}
		\centering 
		\caption{$G=P\square P$, $c=5$}
		\label{tab:5-coloring}
		\begin{tabular}{cccc}
			\toprule
			$i=|H|$ & $n(i)$ & $i=|H|$ & $n(i)$ \\
			\midrule
			4 & 3 & 23 & 12625 \\
			5 & 10 & 24 & 5425 \\
			6 & 22 & 25 & 2666 \\
			7 & 77 & 26 & 6692 \\
			8 & 146 & 27 & 5517 \\
			9 & 238 & 28 & 2956 \\
			10 & 730 & 29 & 988 \\
			11 & 1279 & 30 & 489 \\
			12 & 1627 & 31 & 1139 \\
			13 & 4619 & 32 & 773 \\
			14 & 6691 & 33 & 396 \\
			15 & 7405 & 34 & 127 \\
			16 & 6881 & 35 & 37 \\
			17 & 17196 & 36 & 14 \\
			18 & 19680 & 37 & 26 \\
			19 & 14702 & 38 & 25 \\
			20 & 8497 & 39 & 7 \\
			21 & 21241 & 40 & 2 \\
			22 & 20825 & 41 & 0 \\
			\bottomrule
		\end{tabular}
	\end{minipage}
	\hfill
	\begin{minipage}{0.48\textwidth}
		\centering 
		\caption{$G=P\boxtimes P$, $c=8$}
		\label{tab:8-coloring}
		\begin{tabular}{cccc}
			\toprule
			$i=|H|$ & $n(i)$ & $i=|H|$ & $n(i)$ \\
			\midrule
			4 & 1 & 13 & 90144 \\
			5 & 6 & 14 & 183384 \\
			6 & 24 & 15 & 102816 \\
			7 & 136 & 16 & 127512 \\
			8 & 456 & 17 & 250104 \\
			9 & 1860 & 18 & 35280 \\
			10 & 8064 & 19 & 3144 \\
			11 & 16392 & 20 & 0 \\
			12 & 32568 & &  \\
			\bottomrule
		\end{tabular}
		\vspace{1em} 
		\caption{$G=T_3$, $c=8$}
		\label{tab:8-coloring-T3}
		\begin{tabular}{cccc}
			\toprule
			$i=|H|$ & $n(i)$ & $i=|H|$ & $n(i)$ \\
			\midrule
			4 & 2 & 14 & 387984 \\
			5 & 11 & 15 & 320736 \\
			6 & 44 & 16 & 557112 \\
			7 & 216 & 17 & 1103904 \\
			8 & 756 & 18 & 245520 \\
			9 & 3000 & 19 & 38304 \\
			10 & 13284 & 20 & 1800 \\
			11 & 28872 & 21 & 3480 \\
			12 & 59868 & 22 & 960 \\
			13 & 177384 & 23 & 0 \\
			\bottomrule
		\end{tabular}
	\end{minipage}
\end{table}

The data presented in the tables above support the correctness of Theorems \ref{1} and \ref{2}. The relevant Python code can be downloaded from the open-source project available at 

\url{https://github.com/WentaoZhangT-Zero/Nonrepetitive}.

\section{Discussions}

In the algorithms of this paper, the computer’s running time grows exponentially as two factors increase: the graph’s average degree and the initial number of colors. Additionally, much of this time is spent on data storage and retrieval. Even if we adjust data formats or optimize storage and retrieval methods, it remains difficult to obtain better lower bounds using the current algorithms.

For more general graph classes, getting a good lower bound for the non-repetitive chromatic number with the method in this paper depends on two key choices: which initial graph \(H\) to select, and what order to follow when adding new vertices.

In fact, we also care about the non-repetitive chromatic number of dense graphs—for example, \(K_n \square K_n\). Currently, we only know that for \(\pi(K_n \square K_n)\), its lower bound is of the order of \(n\) and its upper bound is of the order of \(n^2\); these bounds are referenced in Open Problem 3.28 of \cite{WoodDavidR2021NGC} and Question 5.3 of \cite{tao2024nonrepetitivecoloringsgrids}. Algorithms similar to those in this paper may work for cases where \(n\) is small, and could thus help derive a more reasonable conjecture about \(\pi(K_n \square K_n)\).

\section*{Acknowledgment}

We thank Andre Kündgen and Hehui Wu for their helpful discussions and support. Their contributions were greatly appreciated.

\bibliographystyle{plain}
\bibliography{references}

\begin{thebibliography}{10}

\bibitem{alon2002nonrepetitive}
N.~Alon, J.~Grytczuk, M.~Ha{\l}uszczak, and O.~Riordan.
\newblock Nonrepetitive colorings of graphs.
\newblock {\em Random Structures Algorithms}, 21(3-4):336--346, 2002.

\bibitem{barat2007square}
J.~Bar{\'a}t and P.~Varj{\'u}.
\newblock On square-free vertex colorings of graphs.
\newblock {\em Studia Sci. Math. Hungar.}, 44(3):411--422, 2007.

\bibitem{barat2005NoNG}
J.~Bar{\'a}t and D.~R. Wood.
\newblock Notes on nonrepetitive graph colouring.
\newblock {\em Electron. J. Combin.}, 15:R99, 2008.

\bibitem{BOTLER2025221}
F.~Botler, W.~Lomenha, and J.~P. {de Souza}.
\newblock On nonrepetitive colorings of paths and cycles.
\newblock {\em Discrete Appl. Math.}, 360:221--228, 2025.

\bibitem{brevsar2007nonrepetitive}
B.~Bre{\v{s}}ar, J.~Grytczuk, S.~Klav{\v{z}}ar, S.~Niwczyk, and I.~Peterin.
\newblock Nonrepetitive colorings of trees.
\newblock {\em Discrete Math.}, 307(2):163--172, 2007.

\bibitem{currie2002there}
J.~D. Currie.
\newblock There are ternary circular square-free words of length $n$ for $n\ge
  18$.
\newblock {\em Electron. J. Combin.}, 9(1), 2002.

\bibitem{dujmovicVida2020Pghb}
V.~Dujmovi{\'c}, L.~Esperet, G.~Joret, B.~Walczak, and D.~Wood.
\newblock Planar graphs have bounded nonrepetitive chromatic number.
\newblock {\em Adv. Comb.}, 5, 2020.

\bibitem{dujmovic2020planar}
V.~Dujmovi{\'c}, G.~Joret, P.~Micek, P.~Morin, T.~Ueckerdt, and D.~R. Wood.
\newblock Planar graphs have bounded queue-number.
\newblock {\em J. ACM}, 67(4):1--38, 2020.

\bibitem{FertinGuillaume2004Scog}
G.~Fertin, A.~Raspaud, and B.~Reed.
\newblock Star coloring of graphs.
\newblock {\em J. Graph Theory}, 47(3):163--182, 2004.

\bibitem{grytczuk2007nonrepetitive}
J.~Grytczuk.
\newblock Nonrepetitive colorings of graphs-a survey.
\newblock {\em Int. J. Math. Math. Sci.}, 2007(1):1053--1062, 2007.

\bibitem{harant2012nonrepetitive}
J.~Harant and S.~Jendrol.
\newblock Nonrepetitive vertex colorings of graphs.
\newblock {\em Discrete Math.}, 312(2):374--380, 2012.

\bibitem{KESZEGH2020108}
B.~Keszegh and X.~Zhu.
\newblock A note about online nonrepetitive coloring k-trees.
\newblock {\em Discrete Appl. Math.}, 285:108--112, 2020.

\bibitem{kundgen2008nonrepetitive}
A.~K{\"u}ndgen and M.~J. Pelsmajer.
\newblock Nonrepetitive colorings of graphs of bounded tree-width.
\newblock {\em Discrete Math.}, 308(19):4473--4478, 2008.

\bibitem{MARX200913}
D.~Marx and M.~Schaefer.
\newblock The complexity of nonrepetitive coloring.
\newblock {\em Discrete Appl. Math.}, 157(1):13--18, 2009.

\bibitem{tao2024nonrepetitivecoloringsgrids}
T.~Tao.
\newblock The nonrepetitive coloring of grids.
\newblock {\em Discrete Math.}, 349(2):114828, 2026.

\bibitem{thue1906uber}
A.~Thue.
\newblock Uber unendliche zeichenreihen.
\newblock {\em Norske Vid Selsk. Skr. I Mat-Nat Kl.(Christiana)}, 7:1--22,
  1906.

\bibitem{toole}
A.~Toole.
\newblock Repetition-free vertex colorings of grid graphs.
\newblock Master's thesis, California State University San Marcos, 2013.

\bibitem{WoodDavidR2021NGC}
D.~R. Wood.
\newblock Nonrepetitive graph colouring.
\newblock {\em Electron. J. Combin.}, 1000:DS24, 2021.

\end{thebibliography}

\end{document}